\documentclass[12pt]{article}
\usepackage{amssymb}
\newtheorem{theorem}{Theorem}

\newtheorem{proposition}[theorem]{Proposition}
\newtheorem{lemma}[theorem]{Lemma}
\newtheorem{sublemma}[theorem]{Sublemma}
\newtheorem{corollary}[theorem]{Corollary}
\newtheorem{remark}[theorem]{Remark}
\newtheorem{remarks}[theorem]{Remarks}
\newtheorem{example}[theorem]{Example}

\newcommand{\R}{\mathbb{R}}
\newcommand{\Q}{\mathbb{Q}}

\newcommand{\Sf}{\mathbb{S}}

\newcommand{\Les}{\mathbb{L}}

\newcommand{\C}{\mathbb{C}}
\newcommand{\Hy}{\mathbb{H}}

\newcommand{\spa}{\mbox{span}}
\newcommand{\hess}{\mbox{Hess\,}}

\newcommand{\kerl}{\mbox{ker }}

\newcommand{\Py}{\mathbb{P}}

\newcommand{\po}{{\hspace*{-1ex}}{\bf .  }}

\def\Jes{{\cal J}}

\def\<{\langle}

\def\>{\rangle}

\def\va{\varphi}

\def\d{\partial}
\def\bea{\begin{eqnarray*} }
\def\eea{\end{eqnarray*} }
\def\be{\begin{equation} }
\def\ee{\end{equation} }

\def\proof{\noindent{\it Proof: }}
\def\qed{\ifhmode\unskip\nobreak\fi\ifmmode\ifinner
\else\hskip5 pt \fi\fi\hbox{\hskip5 pt \vrule width4 pt
height6 pt  depth1.5 pt \hskip 1pt }}
\begin{document}

\title{All superconformal surfaces in $\R^4$ in \\ terms of   minimal surfaces.
}
\author{Marcos Dajczer and Ruy Tojeiro}
\date{}
\maketitle
\begin{abstract}
We give an explicit construction of any simply-connected superconformal surface
$\phi\colon\,M^2\to \R^4$ in Euclidean space
in terms of a  pair of conjugate minimal surfaces $g,h\colon\,M^2\to\R^4$.
That $\phi$ is superconformal means that its
ellipse of curvature  is a circle at any point. We characterize the pairs $(g,h)$ of
conjugate minimal surfaces that give rise to images of holomorphic curves by an inversion
in $\R^4$ and to images of superminimal surfaces in either a sphere $\Sf^4$ or a hyperbolic
space $\Hy^4$ by an stereographic projection. We also determine the relation
between the pairs $(g,h)$ of conjugate minimal surfaces associated to a superconformal surface
and its image by an inversion. In particular, this yields a new transformation for minimal
surfaces in $\R^4$.
\end{abstract}

\section[Introduction]{Introduction}

For any surface  $\phi\colon\,M^2\to \R^4$ in Euclidean $4$-dimensional space the pointwise
inequality
\be\label{gr}
K+|K_N|\leq\|H\|^2
\ee
holds at every point of $M^2$ \cite{w}. Here $K$ denotes the Gaussian curvature of $M^2$ and
$K_N$ and $H$ are the normal curvature and the mean curvature vector of $\phi$, respectively.
In fact, a similar  inequality was derived by Guadalupe - Rodr\'{\i}guez \cite{gr} for
surfaces of any codimension in space forms of sectional curvature $c$,
namely, $K+|K_N|\leq\|H\|^2+c$. Moreover, it was shown that equality holds
at $p\in M^2$ if and only if the  ellipse of curvature $E(p)$
of $\phi$
at $p$ is a circle.

Recall that the {\it ellipse of curvature\/} of a surface $\phi\colon\,M^2\to N^4$ in a
$4$-dimensional Riemannian manifold at  $p\in M^2$ is the ellipse in the normal space
of $\phi$ at $p$ given by
$$
E(p)=\{\alpha_\phi(X,X)\,:\, X\in T_pM\;\;\mbox{and}\;\; \|X\|=1\},
$$
where $\alpha_\phi$ is the second fundamental form of $\phi$ with
values in the normal bundle; see \cite{gr} and references therein for several interesting
facts on this concept whose study goes back almost a century to the work of  Moore and
Wilson \cite{mw}, \cite{mw2}. Observe that the property of $E(p)$ being a circle is invariant
under conformal changes of the metric of $N^4$.

Following the terminology in \cite{bflpp} we call a surface  $\phi\colon\,M^2\to \R^4$
{\it superconformal\/}  if at any point the ellipse of curvature is a circle.
Thus, the class of superconformal surfaces is invariant under
Moebius transformations of $\R^4$. The condition of superconformality shows up
in several interesting geometric situations. For instance, for a compact oriented surface
integration of (\ref{gr}) over $M^2$ yields the  lower bound
$\int_{M} \|H\|^2\geq 2\pi{\cal X}(M)+|{\cal X}(T^\perp M)|$
due to Wintgen \cite{w} for the Willmore functional  of $\phi$  in terms of the
Euler characteristics of the tangent and normal
bundles. Moreover,  we have equality if and only if $K_N$ does not change sign
and the surface is superconformal.

In this paper,  we provide an explicit construction of any simply connected superconformal
surface in $\R^4$ that is free of minimal and umbilical points.
Start with a simply connected minimal surface
$g\colon\,M^2\to \R^{4}$, oriented by a global conformal diffeomorphism onto either
the complex plane or the unit disk. Then,  consider its conjugate  minimal surface
$h\colon\,M^2\to
\R^{4}$, each of whose components with respect to this global
parameter is the harmonic conjugate of the corresponding component of $g$
(see  \cite{ho} or \cite{la}).  Equivalently, $h_*=g_*\circ J$,
where $J$ is the complex structure on $M^2$ compatible with its orientation.
Notice that  $h$ is determined by $g$  up to a  vector $v\in\R^{\,4}$.

Now, let $\hat{J}_+$ and  $\hat{J}_-$
be the two possible complex structures on $T_g^\perp M$,
 and consider the complex structures $\Jes_+$ and $\Jes_-$ on $g^*T\R^{4}$ given by
$$
\Jes_\pm\circ g_*=g_*\circ J\;\;\;\;\;
\mbox{and}\;\;\;\;\;\Jes_\pm|_{T_g^\perp M}= \hat{J}_\pm.
$$
Our main result reads as follows.

\begin{theorem}\po\label{thm:surf}
Each of the maps $\phi_+\colon\,M^2\to\R^4$ and
$\phi_-\colon\,M^2\to\R^4$  defined by
\be\label{eq:phi3}
\phi_\pm=g+\Jes_\pm h
\ee
parameterizes, at regular points, a superconformal surface.
Moreover, $\phi_+$ and $\phi_-$ are conformal to $g$ and  envelop a  common central
sphere congruence that has
$g$ as its surface of centers. Conversely, any simply connected superconformal
surface that is free of minimal and umbilical points can be constructed in this way.
\end{theorem}

By combining the preceding result with the generalized Weierstrass
parameterization of Euclidean minimal surfaces  \cite{ho},\,\cite{la} we have a parametric
representation of
all simply connected superconformal surfaces in $\R^4$.

Recall that the {\it central sphere congruence\/} (or {\it mean curvature sphere congruence\/})
 of a surface in Euclidean space is the
family of two-dimensional spheres that are tangent to the surface and have the same mean
curvature vector as the surface at the point of tangency. Therefore, our result
implies the fact already known by Rouxel \cite{ro} that superconformal
surfaces $\phi\colon\,M^2\to \R^4$ always arise in pairs $(\phi_+, \phi_-)$ of {\it dual\/}
surfaces that induce conformal metrics on $M^2$ and  envelop a common sphere congruence,
namely, their common central sphere congruence. Hence,  the pair $(\phi_+, \phi_-)$ provides
a solution to the higher codimensional version studied  by Ma \cite{ma1} of the  problem,
first
considered  by Blaschke for surfaces in $\R^3$, of finding all such pairs of  surfaces;
see \cite{hj} for
details on the latter and related facts.

We observe that superconformal minimal surfaces in $\R^4$  are holomorphic curves.
Here and elsewhere, by a surface
$f\colon\,M^2\to\R^{4}$  being holomorphic we mean
that $f$ is complex with respect to a suitable complex
structure of $\R^4$. Therefore,  one obvious way to produce  examples of nonminimal
superconformal  surfaces  is to take compositions of  holomorphic curves with an inversion in
$\R^4$. Notice also that an isolated minimal point can
always be removed by an inversion. Hence, locally and from the point of view of conformal
geometry, assuming that the surface is free of  minimal points in Theorem \ref{thm:surf}
is not essential.

Finally, we point out that  there is a correspondence between holomorphic curves in
$\C\Py^3$ and superconformal surfaces in $\R^4$. In fact,
Theorem $5$ in \cite{bflpp} states that
$\phi\colon\,M^2\to \R^4=\Hy$ is superconformal if and only if either $[\phi,1]$
or $[\bar{\phi},1]$ is the twistor projection under Penrose twistor  fibration
$\pi\colon\,\C\Py^3\to \Hy\Py^1=\Sf^4$ of a holomorphic curve in $\C\Py^3$.
Here $\Hy$ is the set of quaternions, $\Hy\Py^1$ the quaternionic
projective space, $(x,y)\in \Hy^2\setminus (0,0)\mapsto [x,y]\in\Hy\Py^1$
the canonical projection
and $x\in \Hy\mapsto \bar{x}\in \Hy$ the conjugation in $\Hy$.

\vspace{1ex}

In the next result we determine  how the  holomorphic
representative (see  \cite{ho} or \cite{la})
$$
G:=g+ih\colon\,M^2\to\C^4\approx\R^4+i\R^4
$$
of  the minimal surface
$g\colon\,M^2\to \R^4$ associated to an oriented superconformal surface $\phi$ is related
to the holomorphic representative \mbox{$\tilde
G=\tilde{g}+i\tilde{h}$} of  the minimal surface
$\tilde{g}\colon\,M^2\to \R^4$  associated to its composition
$\tilde{\phi}={\cal I}\circ \phi$   with an
inversion ${\cal I}$ in $\R^{4}$ with respect to a sphere
of radius $R$ taken, for simplicity, centered at the origin.

\begin{theorem}\po\label{prop:pairs}
Assume that $\phi$ is not the composition of a holomorphic curve with an inversion. Then
$\tilde G= R^{\,2}T \circ G$, where $T\colon\,\C^4\to \C^4$ is the
holomorphic map
\be\label{eq:T}
T(Z)
=\frac{Z}{\<\!\< Z,Z\>\!\>}
\ee
and $\<\!\< \,\,,\,\,\>\!\>\colon\,\C^{m}\times \C^{m}\to \C$  is the linear inner product on
$\C^{m}$.
\end{theorem}

  For arbitrary $m$, the holomorphic map $T_R:=R^{\,2}T\colon\,\C^{m}\to \C^{m}$ can
  be regarded as the {\em inversion\/} in  $\C^{m}$ with respect to the quadric
$$
\Q_k=\{Z\in\C^{m}: \<\!\< Z,Z\>\!\>=k\}
$$
with $k=R^2$.  Notice that  $T_R$ is defined on $\C^{m}\setminus \Q_0$ and  takes any
quadric $\Q_k$ diffeomorphically onto $\Q_{R^4/k}$.

  As a byproduct of Theorem \ref{prop:pairs}, we obtain the following  remarkable
  property of the map
$T$ that, in particular, yields a new transformation for minimal surfaces in $\R^4$.

\begin{corollary}\po\label{cor:pairs}
 The holomorphic inversion map
$T$ preserves the class of holomorphic curves
$G=g+ih\colon\,M^2\to \C^4$ whose real and imaginary parts $g$
and $h$ define  conjugate minimal immersions into $\R^4$.
\end{corollary}

  Next we characterize the minimal surfaces  that give rise, by means of the construction
  of Theorem~\ref{thm:surf}, to the superconformal surfaces that are images of
holomorphic curves by an inversion in $\R^4$.

\begin{theorem}\po \label{cor:surf2} Let $\phi\colon\,M^2\to \R^4$ be a superconformal
surface parametrized by $(\ref{eq:phi3})$. Then, the  following assertions are equivalent:
\begin{itemize}
\item[$(i)$] The surface $\phi$  is the composition of a
holomorphic curve with an inversion in $\R^4$.
\item[$(ii)$] The  superconformal surface dual to $\phi$ degenerates to a
constant map.
\item[$(iii)$] The  minimal surface $g$  is a holomorphic curve in
$\R^4$.
\end{itemize}
\end{theorem}

We point out that Rouxel \cite{ro} already observed that all spheres of the central sphere
congruence of a composition of a holomorphic curve with an inversion pass through a fixed
point and that, in this case, the surface of centers  is a holomorphic curve.
This is essentially the fact that the first assertion in Theorem \ref{cor:surf2} implies
the remaining two.
\medskip

In the process of proving Theorem \ref{cor:surf2}, the following interesting duality
between holomorphic and anti-holomorphic curves $f\colon\,L^2\to \C^2$  was revealed.
In the next statement
we denote by $f^N$  the normal component of the position vector $f$ in $\R^4\approx\C^2$,
and by ${\cal H}_+$ and ${\cal H}_{-}$ the sets of holomorphic and anti-holomorphic surfaces
in $\C^2$, respectively.

\begin{theorem}\po\label{duality} The map between ${\cal H}_+$ and ${\cal H}_-$ given by
$
f\mapsto f^*=f^N/2\|f^N\|^2
$
is a bijection  such that $(f^*)^*=f$. Moreover, the metrics induced by $f$ and $f^*$ are
conformal.
\end{theorem}

Another class of superconformal surfaces in $\R^4$
is that of stereographic projections of superminimal surfaces in
the sphere $\Sf^4$. Recall that a surface
$g\colon\,L^2\to \Sf^4$ is {\em superminimal\/} if it is minimal
and superconformal; see \cite{fr}
and the references therein for several characterizations of this
very interesting class of surfaces. From a global point of view, it is worth mentioning that
any  minimal immersion of the sphere $\Sf^2$ into $\Sf^4$ is automatically
superminimal \cite{gr}, and that every compact Riemann surface
admits a conformal superminimal immersion into $\Sf^4$ \cite{br}.
The latter result shows that there exist compact
superconformal surfaces in $\R^4$ with arbitrary topology.

A further source of superconformal surfaces
 in $\R^4$ arises by taking the stereographic projections onto a ball in $\R^4$ of
 superminimal surfaces in the hyperbolic space $\Hy^4$. The latter were studied in
\cite{fr} where, in particular, it was proved that there exist complete embedded
simply-connected examples
that are not totally geodesic.

Superminimal surfaces in
$\Sf^4$ and $\Hy^4$ can be regarded as the analogues of holomorphic curves in
$\R^4$. It is natural to ask for  a similar characterization to that of
Theorem \ref{cor:surf2} of the minimal surfaces that produce, by applying the
procedure of Theorem~\ref{thm:surf}, stereographic projections of  superminimal surfaces
in a sphere
$\Sf_R^4=\Sf^4(Re_5;R)$ of radius $R$ centered at $Re_5$, or in a hyperbolic space
 $$\Hy_R^4=\Hy^4(-Re_5;R)=\{X\in \Les^5\,:\,\<X+Re_5,X+Re_5\>=-R^2\}.$$
Here we regard $\R^4$ as the
hyperplane through the origin and normal to the unit vector $e_5$ in either $\R^5$ or
Lorentzian space $\Les^5$, and by the stereographic projection of $\Hy_R^4$ onto $\R^4$
we mean the map that assigns to each $P\in \Hy_R^4$
    the point of $\R^4$ where the line through the points  $-2Re_5$ and $P$
   intersects $\R^4$. Our final result is the following.

\begin{theorem}\po\label{cor:surf3} Let $\phi\colon\,M^2\to \R^4$ be a superconformal
surface parameterized by $(\ref{eq:phi3})$.  Then, the following
assertions are equivalent:
\begin{itemize}
\item[(i)] Either  $\phi$ or its dual  is the composition of a
superminimal surface in a sphere $\Sf_R^4$ (resp., $\Hy_R^4$) with a
stereographic projection of $\Sf_R^4$ (resp., $\Hy_R^4$) onto  $B(0;2R)\subset \R^4$.
\item[(ii)] Both $\phi$ and its dual   are compositions of
superminimal surfaces in a sphere $\Sf_R^4$ (resp., $\Hy_R^4$) with a
stereographic projection of $\Sf_R^4$ (resp., $\Hy_R^4$) onto  $B(0;2R)\subset \R^4$.
\item[(iii)] The holomorphic representative
$G$ of $g$ takes values in  $\Q_{4R^2}$ (resp., $\Q_{-4R^2}$).
\end{itemize}
\end{theorem}

In order to relate the preceding result to Theorem \ref{cor:surf2}, we observe that
holomorphic curves in $\R^4$ can be characterized as the minimal surfaces whose
holomorphic representatives
take values in the quadric $\Q_0$ (see Proposition \ref{le:hol} below).

\section[The ellipse of curvature]{The ellipse of curvature}

In this section we recall some of the basic properties of
the ellipse of curvature of a surface $\phi\colon\,M^2\to N^4$ in a four-dimensional
Riemannian manifold.\vspace{1ex}

 Given an orthonormal basis $\{Y_1,Y_2\}$ of
the tangent space $T_pM$  at  $p\in M^2$,
denote $\alpha_{ij}=\alpha_\phi(Y_i,Y_j)$ for
$1\leq i,j\leq 2$.
 Then, we have for any  $v=\cos \theta
Y_1+\sin\theta Y_2$ that
\be\label{eq:avv}
\alpha_\phi(v,v)=H+\cos2\theta\frac{\alpha_{11}-\alpha_{22}}{2}+\sin2\theta\alpha_{12},
\ee
where $H=\frac{1}{2}(\alpha_{11}+\alpha_{22})$ is the mean
curvature vector of $\phi$ at $p$. This shows that when $v$ goes
once around the unit tangent circle, the vector $\alpha_\phi(v,v)$ goes twice around an
ellipse centered at $H$, the ellipse of curvature $E(p)$ of $\phi$ at $p$.
Clearly, $E(p)$ can degenerate into a line segment or a point.

  It follows from (\ref{eq:avv}) that $E(p)$ is a circle if and only if for some
  (and hence for any) orthonormal basis  of
$T_{\phi(p)}M$ it holds that
\be\label{eq:circ}\<\alpha_{12},\alpha_{11}-\alpha_{22}\>=0\,\,\,\,\,\mbox{and}\,\,\,\,\|\alpha_{11}-\alpha_{22}\|
=2\|\alpha_{12}\|.
\ee

 Let $\{\eta, \zeta\}$ be an orthonormal basis of the normal space $T_{\phi(p)}^\perp M$ at $p$
 with $\eta=H/\|H\|$.
 Take $\{Y_1,Y_2\}$ as an orthonormal tangent
basis   of eigenvectors of the shape operator $A_\zeta$
with respect to $\zeta$, and let $\mu$ and $-\mu$ be the
corresponding eigenvalues. Denoting $\lambda_{ij}=\<A_\eta Y_i,Y_j\>$ for $1\leq i,j\leq 2$
we have
$$\alpha_{11}=\mu\zeta+\lambda_{11}\eta,\,\,\,\,\,\alpha_{12}=\lambda_{12}\eta\,\,\,\,\,\,\,\mbox{and}
\,\,\,\,\,\,\alpha_{22}=-\mu\zeta+\lambda_{22}\eta.$$
Hence, condition (\ref{eq:circ}) for $E(p)$ to be a circle is that $\lambda_{11}=\lambda_{22}$
and $\lambda_{12}=\mu$.
Summarizing, $E(p)$ is a circle if and only if the shape operators $A_\eta$ and $A_{\zeta}$
have the form
\be\label{sffa} \left\{
\begin{array}{l}
A_\eta Y_1=\lambda Y_1+\mu Y_2\vspace{1.5ex}\\
A_\eta Y_2=\mu  Y_1+ \lambda Y_2
\end{array} \right.\;\;\;\;\mbox{and}\;\;\;\;\;
\left\{ \begin{array}{l}
A_\zeta Y_1=\mu Y_1\vspace{1.5ex}\\
A_\zeta Y_2=-\mu Y_2,
\end{array} \right.
\ee
with $\lambda=\|H\|$.  Notice that in this case $\mu$ is  the radius  of $E(p)$.
In particular, $E(p)$ degenerates
into a point if and only if $\mu=0$, that is, $p$ is an umbilical point. If $N^4$ is
a space form, observe that the normal curvature
$K_N=\<R^\perp(Y_1,Y_2)\zeta, \eta\>$ with respect to the oriented orthonormal
bases $\{Y_1,Y_2\}$ and $\{\eta,\zeta\}$ of
$T_pM$ and $T_{\phi(p)}^\perp M$, respectively, is  $K_N=2\mu^2$.

It also follows from (\ref{eq:circ}) that the property that $E(p)$
be a circle is invariant under conformal changes of the metric of
$N^4$. Therefore, any  surface with circular ellipses of
curvature in a (globally) conformally flat $4$-dimensional Riemannian
manifold (in particular any superminimal surface in $\Sf^4$ or $\Hy^4$) gives rise to a
superconformal surface in $\R^4$.

\section[Superconformal surfaces from minimal surfaces]{Superconformal surfaces from
minimal surfaces}

In the two first subsections of this section we prove Theorem
\ref{thm:surf}, starting with the converse statement. In the last one we characterize dual
superconformal surfaces in $\R^4$ that are constructed  from minimal surfaces in $\R^3$
by the procedure of Theorem \ref{thm:surf}.

\subsection[Proof of the converse of Theorem \ref{thm:surf}]{Proof of the converse of
Theorem \ref{thm:surf}}

 Let $\phi\colon\,M^2\to \R^{4}$ be a simply connected oriented surface that has
 nondegenerate circular ellipses of curvature everywhere and nowhere vanishing mean
 curvature vector $H$.
Let $\{\eta, \zeta\}$ be an orthonormal normal frame with
$\eta=H/\|H\|$. By the discussion in the previous section, there
exists an orthonormal tangent frame $\{Y_1,Y_2\}$ such that the
shape operators $A_\eta$ and $A_{\zeta}$ are given everywhere by
(\ref{sffa}), with $\mu$ nowhere vanishing. Thus, we may assume
that $\mu>0$ everywhere. Hereafter,  we let $M^2$ and
$T_\phi^\perp M$  be oriented by the orthonormal frames
$\{Y_1,Y_2\}$ and $\{\eta, \zeta\}$, respectively. Moreover, we
always denote by $J$ the complex structure on  $M^2$ compatible
with its orientation.

We define $g\colon\,M^2\to \R^{4}$  by \be\label{eq:phi2}
g=\phi+r\eta,
\ee
where $r=1/\lambda$. We write
\be\label{eq:zeta}
\zeta=g_*Z+a\xi,
\ee
where $Z\in TM$,  $\xi$ is a unit normal vector
field to $g$ and  $a\in C^\infty(M)$ is such that
\be\label{eq:znorm}
\|Z\|^2+a^2=1.
\ee
Since  $\ker(A_\eta-\lambda I)\neq 0$ everywhere because  $\mu\neq 0$,
we have that $a$ is nowhere vanishing.
Otherwise $g_*Z$ would be somewhere normal to $\phi$,
which would imply, by
taking tangent components for $X=Z$ in
\be\label{eq:phi*}
\phi_*X=g_*X-\<\nabla r,X\>\eta-r\eta_*X,
\ee
that $Z\in\ker(A_\eta-\lambda I)$, and this is a contradiction. Thus, we may assume
that $a>0$ everywhere. Extend $\xi$ to an orthonormal frame $\{\xi,\delta\}$ of
$T_g^\perp M$. It follows from (\ref{eq:phi*}) that
$0=\<\eta,\phi_*X\>=\<\eta,g_*X\>-\<\nabla r,X\>,$
hence the tangent component to $g$ of $\eta$ is $g_*\nabla r$. Thus, we may write
\be\label{eq:eta0}
\eta=g_*\nabla r+\rho\xi+b \delta
\ee
for  $\rho, b\in C^\infty(M)$ satisfying
\be\label{eq:etanorm}
\|\nabla r\|^2+\rho^2+b^2=1.
\ee
\begin{lemma}\po\label{le:basic}  The following holds:
\begin{itemize}
\item[$(i)$] $\rho=0$, $b=\pm a$ and  $JZ=\nabla r$.
\item[$(ii)$] $h=-r\zeta$ satisfies $h_*=g_*\circ J$.
\end{itemize}
\end{lemma}

Before proving Lemma \ref{le:basic}, let us see how it yields  the converse statement
of the theorem. It follows from part $(ii)$ that $g$ and $h$ are conjugate minimal surfaces.
By (\ref{eq:phi2}), (\ref{eq:eta0}) and part~$(i)$, we have \be\label{1}
\phi=g-r\eta=g-rg_*\nabla r\pm ar\delta. \ee Now, \be\label{2}-rg_*\nabla
r=-rg_*JZ=\Jes_\pm
(-rg_*Z)=\Jes_\pm h^T,\ee where $h^T$ denotes the tangent
component to $g$ of the position vector $h$. On the other hand, if
$\hat{J}_+$ and $\hat{J}_-$ are the complex structures on $T_g^\perp M$
defined by
$\hat{J}_\pm\xi=\mp \delta$, then \be\label{3}\pm ar\delta=\hat{J}_\pm(-
ar\xi)=\hat{J}_\pm h^N=\Jes_\pm h^N,\ee where $h^N$ is the normal component
to $g$ of the position vector $h$. We obtain from (\ref{1}), (\ref{2}) and
(\ref{3}) that $\phi$ is given by (\ref{eq:phi3}).
\vspace{1.5ex}

\noindent {\em Proof of Lemma \ref{le:basic}:}
The proof of Lemma \ref{le:basic} will be given in several steps. We start with the following
preliminary facts, where $B_\nu$ denotes the shape operator of $g$ for $\nu\in T_g^ \perp M$.
\begin{sublemma}\po \label{sub1}
  We have
\be\label{bw20} \<B_\delta Z,X\>=a\<\nabla_X^\perp
\delta,\xi\>\,\,\,\,\mbox{for any $X\in TM$} \ee and \be\label{eq:hess} \hess
r(Z)-\frac{1}{r}Z +B_\xi(a\nabla r-\rho Z)+a\nabla \rho=0. \ee
\end{sublemma}
The Codazzi equations for $\phi$ yields
$$
\!\!\!\left\{ \begin{array}{l}
Y_1(\mu)=-\lambda \psi(Y_1)
+\mu\psi(Y_2)+2\mu\Gamma_{22}^1\vspace{1.5ex}\\
Y_2(\mu)=\lambda \psi(Y_2)-\mu\psi(Y_1)
+2\mu\Gamma_{11}^2
\end{array} \right.
\;\;\;\mbox{and}\;\;\;\;
\left\{ \begin{array}{l}
Y_1(\mu)=Y_2(\lambda)+\mu\psi(Y_2)+2\mu\Gamma_{22}^1
\vspace{1.5ex}\\
Y_2(\mu)=Y_1(\lambda)-\mu\psi(Y_1)+2\mu\Gamma_{11}^2,
\end{array} \right.
$$
where $\Gamma_{ii}^j=\<\nabla_{Y_i}Y_i,Y_j\>$ and $\psi$ is the  normal connection
form of $\phi$  given by
$\psi(X)=\<\nabla^\perp_{X}\zeta, \eta\>$.
Subtracting each equation of the first system
from the corresponding equation of the second yields
\be\label{cod3}
\left\{ \begin{array}{l}
Y_2(\lambda)+\lambda\psi(Y_1)=0\vspace{1.5ex}\\
Y_1(\lambda)-\lambda\psi(Y_2)=0. \end{array} \right.
\ee
Differentiating ${\displaystyle g=\phi+\frac{1}{\lambda}\eta}$ and
using (\ref{sffa})
and (\ref{cod3}) gives
$$
\left\{ \begin{array}{l}
{\displaystyle g_*Y_1
=-\frac{1}{\lambda}\left(\psi(Y_2)\eta
+\mu \phi_*Y_2+\psi(Y_1)\zeta\right)
}\vspace{1.5ex}\\
{\displaystyle g_*Y_2=
-\frac{1}{\lambda}\left(-\psi(Y_1)\eta
+\mu\phi_*Y_1+\psi(Y_2)\zeta\right).
}\end{array} \right.
$$
Therefore,  the vector field $
\psi(Y_1)\phi_*Y_1-\psi(Y_2)\phi_*Y_2+\mu\eta $ is normal to $g$,
and since it is orthogonal to $\xi$, it is in the direction of $\delta$.
Thus, from
$$
\left\{ \begin{array}{l}
\zeta_*Y_1=-\mu \phi_*Y_1+\psi(Y_1)\eta
\vspace{1.5ex}\\
\zeta_*Y_2=\mu \phi_*Y_2+\psi(Y_2)\eta,
\end{array} \right.
$$
we obtain that
$\<\zeta_*Y_j,\delta\>=0$,
which is easily seen to be equivalent to (\ref{bw20}). Then,
\be\label{eq:zeta*}
\zeta_*X=g_*DX+\<K,X\>\xi,
\ee
where
\be\label{eq:dx}
DX=\nabla_X Z-aB_\xi
X\;\;\;\;\mbox{and}\;\;\;\;\;K=\nabla a+B_\xi Z.
\ee
The orthogonality
between $\eta$ and $\zeta$ yields
\be\label{eq:zr}
\<Z,\nabla r\>+a\rho=0.
\ee Hence,
\begin{eqnarray}\label{eq:nzx}
\<\nabla_XZ,\nabla r\>\!\!&=&\!\!X\<Z,\nabla r\>-\<Z,\hess r(X)\>\nonumber\\
\!\!&=&\!\!-X(a)\rho -aX(\rho)-\<Z,\hess r(X)\>\nonumber\\
\!\!&=&\!\!-\<\rho\nabla a+a\nabla \rho+\hess r(Z),X\>.
\end{eqnarray}
It follows from (\ref{eq:zeta*}),  (\ref{eq:dx}) and (\ref{eq:nzx}) that \be\label{eq:z1}
\<\zeta_*X,\eta\>=\<DX,\nabla r\>+\rho\<K,X\>=-\<\hess r(Z)+B_\xi(a\nabla
r-\rho Z)+a\nabla \rho,X\>.\ee On the other hand,
\be\label{eq:z2}\<\zeta,\phi_*X\>=X\<\zeta,g\>-\<\zeta_*X,g-r\eta\>=\<Z,X\>
+r\<\zeta_*X,\eta\>,\ee and thus (\ref{eq:hess}) follows from (\ref{eq:z1}), (\ref{eq:z2})
and the fact that $\zeta$ is
normal to $\phi$.\vspace{1ex}

   The next step is to express (\ref{sffa}) in terms of $g$. It is
convenient to use the
orthonormal frame
$$
X_1=\frac{1}{\sqrt{2}}(Y_1+Y_2),\;\;\;\;\;\; X_2=\frac{1}{\sqrt{2}}(Y_1-Y_2),
$$
with respect to which  (\ref{sffa}) becomes
$$
\left\{ \begin{array}{l}
A_\eta X_1=(\lambda+\mu) X_1\vspace{1.5ex}\\
A_\eta X_2=(\lambda-\mu) X_2
\end{array} \right. \;\;\;\;\mbox{and}\;\;\;\;\;
\left\{ \begin{array}{l}
A_\zeta X_1=\mu X_2\vspace{1.5ex}\\
A_\zeta X_2=\mu X_1.
\end{array} \right.
$$
Hence,
\be\label{sff}
\left\{ \begin{array}{l}
\eta_*X_1=-(\lambda
+\mu )\phi_*X_1-\psi(X_1)\zeta
\vspace{1.5ex}\\
\eta_*X_2=-(\lambda
-\mu )\phi_*X_2-\psi(X_2)\zeta
\vspace{1.5ex}\\
\zeta_*X_i=-\mu \phi_*X_j+\psi(X_i)\eta,
\;\; 1\leq i\neq j\leq 2,\\
\end{array} \right.
\ee
In view of (\ref{eq:phi*}), this is equivalent to \be\label{sff2}
\left\{
\begin{array}{l}r^2\mu\eta_*X_1 =\theta_1(g_*X_1-r_1\eta)+r\psi(X_1)\zeta
\vspace{1.5ex}\\
r^2\mu\eta_*X_2
=\theta_2(-g_*X_2+r_2\eta)-r\psi(X_2)\zeta
\vspace{1.5ex}\\
\zeta_*X_i=-\mu g_*X_j+\mu( r\eta_*X_j +r_j\eta)+\psi(X_i)\eta,\;\;
1\leq i\neq j\leq 2,\\
\end{array} \right.
\ee
where $\theta_1=(1+r\mu)$,  $\theta_2=(1-r\mu)$ and $r_i=\<\nabla r,
X_i\>$ for $1\leq i\leq 2$.

We have \be\label{eq:eta*} \eta_*X=g_*QX
+\<T,X\>\xi+\<P,X\>\delta \ee where
$$
\left\{ \begin{array}{l}
{\displaystyle Q=\hess r-\rho B_\xi -b B_\delta}\label{eq:qw}\vspace{1ex}\\
{\displaystyle T=\nabla\rho +B_\xi\nabla
r+\frac{b}{a}B_\delta Z} \label{eq:omega}\vspace{1ex}\\
{\displaystyle P=\nabla b +B_\delta\nabla
r-\frac{\rho}{a}B_\delta Z}. \label{eq:alpha}\end{array} \right.
$$
To proceed we use that the  normal connection form $\psi$ of $\phi$ can be written as
\be\label{eq:ncon} \psi(X)=-\frac{1}{r}\<Z,X\>\,\,\,\,\mbox{for any $X\in TM$} \ee
in terms of  data
related to $g$. This follows immediately from (\ref{eq:hess}) and (\ref{eq:z1}).

Using (\ref{eq:ncon}),  the $\delta$-component of
(\ref{sff2}) gives \be\label{eq:wcomp} \left\{ \begin{array}{l}
r^2\mu\<P,X_1\>=-\theta_1b r_1\vspace{1.5ex}\\
r^2\mu\<P,X_2\>=\theta_2b r_2\vspace{1.5ex}\\
r^2\mu\<P,X_i\>=-r\mu b r_i+b\<Z,X_j\>,
\;\; 1\leq i\neq j\leq 2.\\
  \end{array} \right. \ee
Replacing the first two equations into the last two yields
\be\label{eq:gradrZ}
r_1=-\<Z,X_2\> \;\;\;\mbox{and}\;\;\;\;\;r_2=\<Z,X_1\>.
\ee
Taking the tangent component to $g$ of (\ref{sff2}) and using
(\ref{eq:gradrZ}) we obtain
\be\label{sff3}
\left\{ \begin{array}{l}
r^2\mu QX_1-\theta_1SX_1+r_2Z=0\vspace{1.5ex}\\
r^2\mu QX_2+\theta_2SX_2+r_1Z=0\vspace{1.5ex}\\
rDX_1+r\mu SX_2-r^2\mu QX_2+r_2\nabla r=0\vspace{1.5ex}\\
rDX_2+r\mu SX_1-r^2\mu QX_1-r_1\nabla r=0,\\
\end{array} \right.
\ee
where
\be\label{eq:S}
S=I-\<\nabla r, *\>\nabla r.
\ee
Finally, computing  the $\xi$-component of (\ref{sff2}) yields
\be\label{eq:xicomp} \left\{ \begin{array}{l}
r^2\mu\<T,X_1\>=-\theta_1\rho r_1 -a\<Z,X_1\>\vspace{1.5ex}\\
r^2\mu\<T,X_2\>=\theta_2\rho r_2 +a\<Z,X_2\>\vspace{1.5ex}\\
r\<K,X_i\>=r^2\mu\<T,X_j\>+ r\mu\rho r_j-\rho\<Z,X_i\>,\;\; 1\leq i\neq
j\leq 2.\\
\end{array} \right.
\ee

We now prove:
\begin{sublemma} {\hspace*{-1ex}}\textbf{. } \label{sub2} The metrics
induced by $g$ and
$\phi$ are conformal. Namely, \be\label{eq:conf}
\<g_*X,g_*Y\>=\frac{r^2\mu^2}{a^2}\<\phi_*X,\phi_*Y\>.
\ee
\end{sublemma}
 From (\ref{eq:phi*}) and (\ref{eq:eta*}) we have
\be\label{eq:orth2}
\begin{array}{l}
\delta_{ij}=\<\phi_*X_i,\phi_*X_j\>=\<X_i,X_j\> -r_ir_j-2r\<QX_i,X_j\>
\vspace{1ex}\\\hspace*{10ex} + r^2(\<QX_i,QX_j\>+\<T,X_i\>\<T,X_j\>
+\<P,X_i\>\<P,X_j\>).
\end{array}\ee
Taking inner products of the first and second  equations in
(\ref{sff3}) by $X_2$ and $-X_1$, respectively,
and adding them up taking (\ref{eq:gradrZ}) into account, yields
$$
\<X_1,X_2\>=0.
$$
We compute from the first two equations in (\ref{sff3}), bearing in
mind (\ref{eq:gradrZ}), that
\be\label{eq:qwuiuj}
\left\{\begin{array}{l}
r^2\mu\<QX_1,X_1\>=\theta_1(\|X_1\|^2-r_1^2)-r_2^2
\vspace{1.5ex}\\
r^2\mu\<QX_2,X_2\>=-\theta_2(\|X_2\|^2-r_2^2)+r_1^2
\vspace{1.5ex}\\
r\<QX_1,X_2\>=-r_1r_2.\\
\end{array} \right.
\ee
Using (\ref{eq:znorm}) and (\ref{eq:zr}), we have
\be\label{eq:qwuiqwuj}
\hspace*{-2ex}\left\{\begin{array}{l}
\!\!r^4\mu^2\|QX_1\|^2=\theta_1^2(\|X_1\|^2+(\|\nabla
r\|^2-2)r_1^2)-2\theta_1(r_2^2+a\rho r_1r_2)
+(1-a^2)r_2^2\vspace{1.5ex}\\
\!\!r^4\mu^2\|QX_2\|^2=\theta_2^2(\|X_2\|^2+(\|\nabla
r\|^2-2)r_2^2)-2\theta_2(r_1^2-a\rho r_1r_2)
+(1-a^2)r_1^2\vspace{1.5ex}\\
\!\!r^4\mu^2\<QX_1,QX_2\>=(\theta_1\theta_2
(1+\rho^2+b^2)-\theta_1-\theta_2-a^2+1)r_1r_2
-a\rho(\theta_1r_1^2-\theta_2r_2^2).\\
\end{array} \right.
\ee
  From (\ref{eq:gradrZ}) and (\ref{eq:xicomp}) we obtain
$$ r^2\mu\<T,X_1\>=-\theta_1\rho
r_1-ar_2\;\;\;\mbox{and}\;\;\; r^2\mu\<T,X_2\>=\theta_2\rho r_2-ar_1. $$
Thus, \be\label{eq:twui2} \left\{ \begin{array}{l} {\displaystyle
r^4\mu^2\<T,X_1\>^2=\theta_1^2\rho^2r_1^2+a^2r_2^2 +2\theta_1a\rho
r_1r_2}\vspace{1.5ex}\\{\displaystyle r^4\mu^2\<T,X_1\>\<T,X_2\>=
(a^2-\theta_1\theta_2\rho^2)r_1r_2 +\theta_1a\rho r_1^2}-\theta_2a\rho
r_2^2
\vspace{1.5ex}\\
{\displaystyle
r^4\mu^2\<T,X_2\>^2=\theta_2^2\rho^2r_2^2+a^2r_1^2-2\theta_2a\rho
r_1r_2}.\end{array} \right.
\ee
  From the first two equations in (\ref{eq:wcomp}) we get
\be\label{eq:pwui2}
\left\{ \begin{array}{l}
{\displaystyle r^4\mu^2\<P,X_1\>^2=\theta_1^2b^2r_1^2}
\vspace{1.5ex}\\
{\displaystyle r^4\mu^2\<P,X_2\>^2=\theta_2^2b^2r_2^2}
\vspace{1.5ex}\\
{\displaystyle r^4\mu^2\<P,X_1\>\<P,X_2\>=-\theta_1\theta_2
b^2r_1r_2.}\end{array} \right.
\ee
Replacing (\ref{eq:qwuiuj}), (\ref{eq:qwuiqwuj}), (\ref{eq:twui2}) and
(\ref{eq:pwui2})
into (\ref{eq:orth2})  we end up with
$$
\|X_1\|^2=r^2\mu^2+r_1^2+r_2^2=\|X_2\|^2,
$$
and (\ref{eq:conf}) follows easily.\vspace{1.5ex}

It follows from Sublemma \ref{sub2}
and (\ref{eq:gradrZ})  that $JZ=\nabla r$. We conclude from
(\ref{eq:zr}) that $\rho=0$, hence $b=\pm a$ by
(\ref{eq:znorm}) and (\ref{eq:etanorm}), and the proof of  $(i)$ is completed.\vspace{1ex}\\
We now prove $(ii)$. Replacing the first two equations of (\ref{sff3})
into the last two gives
\be\label{ds}
\left\{ \begin{array}{l}
rDX_1+SX_2+r_1Z+r_2\nabla r=0\vspace{1.5ex}\\
rDX_2-SX_1+r_2Z-r_1\nabla r=0,\\
\end{array} \right.
\ee
that can be written as
\be\label{eq:D}
rDX=-JX-\<\nabla r,X\>Z.
\ee
On the other hand, replacing the first two equations of
(\ref{eq:xicomp}) into the last
two yields
\be\label{eq:rKui}
\left\{ \begin{array}{l}
r\<K,X_1\>=a\<Z,X_2\>
\vspace{1.5ex}\\
r\<K,X_2\>=-a\<Z,X_1\>.\\
\end{array} \right.
\ee
Taking (\ref{eq:gradrZ}) into account, the preceding equations reduce to
\be\label{eq:nablaar}
rB_\xi Z+\nabla(ar)=0.
\ee
 From (\ref{bw20}) we have
$$
\tilde{\nabla}_X\xi=-g_*B_\xi X+\nabla_X^\perp \xi\vspace{1ex}\\
=-g_*B_\xi X-\frac{1}{a}(\alpha_g(Z,X)-\<B_\xi Z,X\>\xi),
$$
where $\alpha_g$ denotes the second fundamental form of $g$. Hence,
$$-ar\tilde{\nabla}_X\xi+r\<B_\xi Z,X\>\xi=arg_*B_\xi X+r\alpha_g(Z,X).$$
In view of (\ref{eq:nablaar}) the left-hand-side is
$\tilde{\nabla}_X(-ar\xi)$.
For the right-hand-side we have
\begin{eqnarray*}
arg_*B_\xi X+r\alpha_g(Z,X)\!\!&=&\!\!arg_*B_\xi
X+r(\tilde{\nabla}_Xg_*Z-g_*\nabla_XZ)\vspace{1ex}\\
\!\!&=&\!\! g_*(arB_\xi X-r\nabla_XZ-X(r)Z)+\tilde{\nabla}_X(rg_*Z).
\end{eqnarray*}
Therefore, we obtain using (\ref{eq:D}) that
$$
h_*X=g_*(arB_\xi X-r\nabla_XZ-X(r)Z)=g_*(-rDX-X(r)Z)=g_*JX.\qed
$$

\subsection[Proof of the direct statement  of Theorem \ref{thm:surf}]
{Proof of the direct statement of Theorem \ref{thm:surf}}

   For the proof of the direct statement  we need the following
   general fact.

\begin{proposition}\po\label{le:conj}
Let $g\colon\,M^2\to \R^{n+2}$  be a simply connected oriented minimal surface with complex
structure $J$
compatible with the orientation and let $h\colon\,M^2\to \R^{n+2}$ be a conjugate
minimal surface such that
 \be\label{conj}h_*=g_*\circ J.
\ee   Then $r=\|h\|$ satisfies
$\|\nabla r\|\leq 1$
everywhere. Moreover, on the complement of the subset of isolated
points of $M^2$ where $a=\sqrt{1-\|\nabla r\|^2}$ vanishes, there
exists a  smooth unit normal vector field $\xi$ to $g$ such that
$$
h=-r(g_*\nabla r+a \xi).
$$
Furthermore,
\be\label{bw201} \<B_\delta J\nabla r,X\>+a\<\nabla_X^\perp
\delta,\xi\>=0\,\,\,\mbox{ for all}\,\,\,\delta\in \spa\{\xi\}^\perp\ee and
\be\label{eq:bxi}
B_\xi=\frac{1}{ar}(r\mbox{Hess\,} r-S)\circ J,
\ee
where $S$ is given by (\ref{eq:S}).\end{proposition}
\proof Decompose  $h=g_*T+h^N$ into tangent and normal components to $g$.
 From (\ref{conj}) we obtain \be\label{eq:combe}
\left\{ \begin{array}{l}
\nabla_XT-B_{h^N} X=JX
\vspace{1.5ex}\\
\alpha_g(X,T)+\nabla_X^\perp h^N=0.
\end{array} \right.
\ee
It also follows from (\ref{conj}) that the tangent components of the position vector $h$
with respect to $g$ and $h$
coincide. Since the latter is $h_*(r\nabla r)$,
we get
$$
g_*T=h_*(r\nabla r)=g_*(rJ\nabla r),
$$
 hence $T=rJ\nabla r$. From
$\|h\|^2=\|T\|^2+\|h^N\|^2$ we get
$
r^2=r^2\|\nabla r\|^2+\|h^N\|^2,
$
which implies that $\|\nabla r\|\leq 1$
holds everywhere.  By the real analyticity of $g$ and $h$ the points
where the function $a$ vanishes are isolated. On the complement of
the subset of such points we have $\|h^N\|=ar$. Thus, we can write
$h^N=-ar\xi$ for a  unit normal vector field $\xi$. Then, using that
$J\circ B_\xi = -B_\xi \circ J$, for $B_\xi$ is traceless, the
first equation in (\ref{eq:combe}) reduces to (\ref{eq:bxi}),
whereas the $\spa\{\xi\}^\perp$-component of the second yields
(\ref{bw201}).
\vspace{1,5ex}\qed

Setting $\delta_\pm=\hat{J}_\pm\xi$ and $\eta_\pm=g_*\nabla r + a\delta_\pm$,
we have from Proposition \ref{le:conj} that
$$
\phi_\pm=g-r\eta_\pm.
$$
It follows from (\ref{eq:phi*})  that $\eta_\pm$ is a unit normal vector field to $\phi_\pm$.
Let  $\zeta$  be defined by (\ref{eq:zeta}) with $Z=-J\nabla r$. Then $\zeta$  has unit
length and is orthogonal
to $\eta_\pm$. We obtain from (\ref{bw201}) that (\ref{eq:zeta*}) holds, hence
we have (\ref{eq:z1})  with $\rho=0$, and  also  (\ref{eq:z2}).
>From (\ref{eq:bxi}) we get
\be\label{hess2}
r\hess r(Z) -Z+arB_\xi\nabla r=0,
\ee
which implies, using (\ref{eq:z1}) (with $\rho=0$)  and (\ref{eq:z2}),
that $\zeta$ is normal to $\phi_\pm$.

 Therefore, to complete the proof it suffices to show that there exists an orthonormal
 tangent frame $\{X_1,X_2\}$
 (with respect to the metric induced by $\phi_\pm$)
satisfying (\ref{sff}). Since $B_{\delta_\pm}$ and $B_\xi$ are traceless symmetric
$2\times 2$ matrices, we have

\be\label{eq:alpha2}(B_{\delta_\pm}+B_\xi J)^2
=\alpha_\pm^2 I
\ee
for some smooth functions $\alpha_\pm$. By analyticity, either $\alpha_\pm$ vanishes
identically or it vanishes only
 at isolated points. In the first case, a standard argument shows that the complex
 structure ${\cal J}=J\oplus
\hat{J}_\pm$ on $g^*T\R^4$ is a parallel tensor,
hence  defines a complex structure on $\R^4$ with respect to which  $g$ is
holomorphic. Then, in this case the
 conclusion follows from Theorem
\ref{cor:surf2}. Therefore, we may assume in the sequel that
$\alpha_\pm$ is nowhere vanishing, hence there exists $\mu_\pm\in C^\infty(M)$
such that $\alpha_\pm=-a/r^2\mu_\pm$. Since
 $B_{\delta_\pm}+B_\xi J=\alpha_\pm R_\pm$ for some reflection
$R_\pm$ by (\ref{eq:alpha2}), it follows using (\ref{eq:bxi}) that  \be\label{eq:bw2}
B_{\delta_\pm}=\frac{1}{a}(\hess r-\frac{1}{r}S)-\frac{a}{r^2\mu_\pm}R_\pm.
\ee

 Let $\{\bar{X}^\pm_1, \bar{X}^\pm_2\}$ be the orthonormal basis of eigenvectors of $R_\pm$
 (with respect to the
metric induced by $g$), with $\bar{X}^\pm_1$ corresponding
to the eigenvalue $+1$ and $\bar{X}^\pm_2=J\bar{X}^\pm_1$. Define
$$
X_j^\pm=\frac{r\mu_\pm}{a}\bar{X}^\pm_j,\;\;\; 1\leq j\leq 2.
$$
 We claim that $\{X_1^\pm,X_2^\pm\}$  is the desired orthonormal frame.
In order to prove (\ref{sff}), it  suffices to show that (\ref{sff2}), or equivalently,
(\ref{eq:wcomp}),
(\ref{sff3}) and (\ref{eq:xicomp}), holds for $X_1^\pm$ and $X_2^\pm$.

  Since we have (\ref{eq:gradrZ}), because  $JX_1^\pm=X_2^\pm$ and $JZ=\nabla r$,
  system (\ref{eq:wcomp})
  reduces to its first two equations. These are in turn equivalent to
$$
rB_{\delta_\pm}\nabla r+\frac{a}{r\mu_\pm}R_\pm\nabla r+\nabla(ar)=0,
$$
which follows from (\ref{eq:bw2}).
Now, (\ref{eq:bw2}) also implies that
$$rQ=S+\frac{a^2}{r\mu_\pm}R_\pm.$$
 Moreover,  from (\ref{eq:bxi}) we get (\ref{eq:D}), hence (\ref{sff3}) is satisfied.

 From (\ref{eq:bxi}) we also obtain (\ref{eq:nablaar}), and hence (\ref{eq:rKui}).
Moreover, (\ref{eq:bxi}) and (\ref{eq:bw2}) imply  that
$$
B_\xi\nabla r+B_{\delta_\pm} Z+\frac{a}{r^2\mu_\pm}R_\pm Z=0,
$$ thus
(\ref{eq:xicomp}) is satisfied.
Since we now have (\ref{eq:qwuiuj}), (\ref{eq:qwuiqwuj}), (\ref{eq:twui2}) and
(\ref{eq:pwui2}), then $$
\<{\phi_\pm}_*X_i^\pm,{\phi_\pm}_*X_j^\pm\>=\delta_{ij}
$$
follows by using that $\<X_1^\pm,X_2^\pm\>=0$ and $\|X_j^\pm\|^2=r^2\mu_\pm^2/{a}^2$
for $1\leq j\leq 2$.

  Finally, that $\phi_+$ and $\phi_-$ envelop a common central sphere congruence,
  with $g$ as the surface of centers, follows from the facts that for each $p\in M^2$
  we have
$$
\phi_+(p)+H_+(p)/\|H_+(p)\|^2
=g(p)=\phi_-(p)+H_-(p)/\|H_-(p)\|^2
$$
and  $$(\phi_+)_*T_pM\oplus \spa\{H_+(p)\}
=(\phi_-)_*T_pM\oplus \spa\{H_-(p)\},
$$
for  $\zeta(p)$ is orthogonal to both subspaces.\qed

\begin{remarks}\po\label{reg} {\em $(1)$ It follows from the proof that  $\phi_\pm$
may fail to be regular only  at  points where the function $a$ in Lemma
\ref{le:conj} vanishes, that is, at points where the position vector of $h$ is tangent
to $g$, and at points where the
shape operator of $g$ with respect to  any normal direction $\beta$  satisfies
$B_{\hat{J}_\pm\beta}=-B_\beta\circ J$.
The latter can be seen as the  ``holomorphic" points of $g$, that is, points where the
ellipse of curvature of $g$ is a
circle.\vspace{1ex}\\
$(2)$ If we change the conjugate minimal
surface $h$ by $h+v$ for any $v\in\R^4$, then the
corresponding surfaces $\phi_\pm$ are changed by addition of
$\Jes_\pm v$. One can check that  the latter is just a
parameterization  of an open subset of the two-dimensional equator
in $\Sf^3(\|v\|)$ orthogonal to $v$. Moreover, the
parameterization is conformal to $g$ and singularities occur at
points where the ellipse of curvature of $g$ is a  circle.\vspace{1ex}\\
$(3)$ It follows from Sublemma \ref{sub2} that the
 metrics $\<\,,\,\>_+$ and $\<\,,\,\>_-$ induced on $M^2$ by $\phi_+$ and $\phi_-$,
 respectively, are related by
$\mu^2_+\<\,,\,\>_+=\mu^2_-\<\,,\,\>_-.$\vspace{1ex}\\
$(4)$ For any element of the associated family
$g_\theta=\cos\theta g+\sin\theta h$
of the minimal surface $g$, we have a pair $(\phi_+^\theta,\phi_-^\theta)$
of dual superconformal surfaces in $\R^4$ determined by the pair
$(g_\theta,h_\theta)$ of conjugate minimal surfaces $g_\theta$ and
$h_\theta=-\sin\theta g+\cos\theta h$, which also satisfy
${h_\theta}_*={g_\theta}_*\circ J$.
Namely,
$$
\phi_\pm^\theta=g_\theta+\Jes_\pm h_\theta=\Jes_\pm^\theta
g+\tilde{\Jes}_\pm^\theta h,
$$
where
$$
\Jes_\pm^\theta\circ g_*=g_*(cos\theta I-\sin\theta
J),\,\,\,\,\,\,\,\,\tilde{\Jes}_\pm^\theta\circ g_*=g_*(\sin\theta I+\cos\theta J)
$$
and ${\Jes_\pm^\theta}|_{T_g^\perp M}$ and
${\tilde{\Jes}_\pm^\theta}|_{T_g^\perp M}$ are given, respectively,  by
$$
{\Jes_\pm^\theta}|_{T_g^\perp M}=cos\theta I-\sin\theta
{\Jes}_\pm,\,\,\,\,\,\,\,\,
{\tilde{\Jes}_\pm^\theta}|_{T_g^\perp M}=\sin\theta I+\cos\theta
{\Jes}_\pm.
$$
It is an interesting question whether $\{\phi_\pm^\theta\}$ coincides with the associated
family of $\phi_\pm$ in the sense of \cite{ma3}, Corollary 2.7.

}\end{remarks}

 \subsection[Superconformal surfaces from minimal surfaces in $\R^3$]{Superconformal
 surfaces from minimal surfaces in $\R^3$}

 In this subsection  we consider the problem of determining when a pair of  dual superconformal
 surfaces $\phi_+$ and $\phi_-$ differ by an inversion in $\R^4$.

\begin{proposition}\po\label{le:conf} Let $(\phi_+,\phi_-)$ be a pair of dual
superconformal surfaces constructed  from a minimal
surface $g\colon\,L^2\to \R^3\subset\R^4$ by the procedure of Theorem \ref{thm:surf}.
Then, the following holds:
\begin{itemize}
\item[(i)] The maps
$\phi_+$ and $\phi_-$ differ by a reflection with respect to $\R^3$.
\item[(ii)] For any inversion ${\cal I}$  with respect to a sphere centered at a point
$P_0\not\in \R^3$, the maps ${\cal I}\circ \phi_+$ and ${\cal I}\circ \phi_-$ define dual
superconformal surfaces that differ by an inversion with respect to the sphere ${\cal I}(\R^3)$.
\end{itemize}
 Conversely, any pair of dual superconformal surfaces that differ by an inversion in
 $\R^4$ arises as in $(ii)$.
\end{proposition}


\proof If $g\colon\,L^2\to \R^4$ is a minimal surface such that $g(L^2)\subset \R^3$, then
$$
\phi_\pm=g+Jh^T\pm \<h,N\>e_4,
$$
where $N$ is a unit normal vector field to $g$ in $\R^3$. This proves the first assertion.
The second follows from a well known property of inversions. For the converse,
if a pair $(\phi_+, \phi_-)$ of dual superconformal surfaces  differ by an inversion
in $\R^4$ with respect to a hypersphere $S$, then each element of their common central
sphere congruence is orthogonal to $S$, since it passes through a pair of inverse points.
Then the image of $S$ by an inversion ${\cal I}$ with respect to a point contained in it
is a hyperplane $\R^3$ and any element of the common central sphere congruence of
${\cal I}\circ \phi_+$ and  ${\cal I}\circ \phi_-$, being orthogonal to $\R^3$,
has its center therein. Therefore ${\cal I}\circ \phi_+$ and  ${\cal I}\circ \phi_-$
are constructed as in Theorem \ref{thm:surf} from a minimal
surface $g\colon\,L^2\to \R^3$.\vspace{1.5ex}\qed

The following example shows that minimal surfaces in $\R^3$ give rise to nontrivial examples
of superconformal surfaces in $\R^4$ by means of the construction in Theorem \ref{thm:surf}.

\begin{example}\label{ex:cathel}\po{\em
Consider the catenoid and the helicoid in $\R^3$ parameterized, respectively, by
$
g(u,v)=(\cosh v\cos u,\cosh v\sin u,v)
$
and
$
h(u,v)= (-\sinh v\sin u, \sinh v\cos u, -u).
$
Then the corresponding dual superconformal surfaces given by Theorem \ref{thm:surf} are
$$
\phi_\pm= \frac{1}{\cosh v}
(\cos u - u\sin u, \sin u +u\cos u,v\cosh v-\sinh v,\pm
u\sinh v).
$$
}
\end{example}

\section[Proof of Theorem \ref{prop:pairs}]{Proof of Theorem \ref{prop:pairs}}

    For the proof of Theorem \ref{prop:pairs} we need the following well-known  fact.

\begin{lemma}\po\label{le:inv}  Let $f\colon\,M^n\to \R^{N}$ be an isometric immersion and
let ${\cal I}$ be an
inversion with respect to a sphere with radius $R$ centered at $P_0\in \R^{N}$. Then,
\be \label{eq:vbi}
{\cal P}\xi=\xi-2\frac{\<f-P_0,\xi\>}{\<f-P_0,f-P_0\>}(f-P_0)
\ee
is a vector bundle isometry between the normal bundles $T_f^\perp M$ and
$T_{{\cal I}\circ f}^\perp M$. Moreover,
the shape operators $A_\xi$ and $\tilde{A}_{{\cal P}\xi}$ are related by
\be\label{eq:sop}
\tilde{A}_{{\cal P}\xi}
=\frac{1}{R^2}\left(\<f-P_0,f-P_0\>A_\xi+2\<f-P_0,\xi\>I\right).
\ee
\end{lemma}

In Section $6$ we will need the following analogue of Lemma \ref{le:inv}  for isometric
immersions $f\colon\,M^n\to \Les^{N}$ into Lorentzian space and the
``inversion"
$$
{\cal I}(P)=P_0-\frac{R^2}{\<P-P_0,P-P_0\>}(P-P_0),\;\;\;
P\neq P_0.
$$
with respect to the hyperbolic space $\Hy^{N-1}_R:=
\Hy^{N-1}(P_0;R)$ of radius $R$ ``centered" at $P_0\in \Les^N$,~i.e.,
$$
\Hy^{N-1}_R=\{P\in \Les^N\,:\,\<P-P_0,P-P_0\>=-R^2\}.
$$

\begin{lemma}\po\label{le:invv}  Let $f\colon\,M^n\to \Les^{N}$ be an
isometric immersion and let ${\cal I}$ be an
inversion with respect to $
\Hy_R^{N-1}$.
Then (\ref{eq:vbi}) and (\ref{eq:sop}) hold true if
we replace $R^2$ by $-R^2$ in the latter.
\end{lemma}

\proof   We provide a proof for the sake of completeness, which also applies for
Lemma \ref{le:inv}.
An easy computation shows that
\be\label{eq:tang}\tilde{f}_*=-\frac{R^2}{\<f-P_0,f-P_0\>}{\cal P}\circ
f_*,\ee
where ${\cal P}\colon\,f^*T\Les^{N}\to f^*T\Les^{N}$ is given by
$$
{\cal P}Z=Z-2\<f-P_0,Z\>\<f-P_0,f-P_0\>^{-2}(f-P_0).
$$
Since ${\cal P}$ is easily seen to be a vector bundle isometry, it
follows from (\ref{eq:tang}) that it restricts to a  vector bundle
isometry of $T_f^\perp M$ onto itself. Denoting by  $\bar{\nabla}$  the derivative
in $\Les^N$,  equation (\ref{eq:sop})
follows by taking tangent components in
$$
\begin{array}{l} -\tilde{f}_*\tilde{A}_{{\cal
P}\xi}X+\tilde{\nabla}^\perp {\cal P}\xi=
{\displaystyle \bar{\nabla}_X\tilde{f}_*Y
=\bar{\nabla}_X(\xi
-2\frac{\<f-P_0,\xi\>}{\<f-P_0,f-P_0\>}(f-P_0))}
\vspace{1.5ex}\\\hspace*{19.7ex}={\displaystyle
-{\cal P}f_*(A_\xi+2\frac{R^2}{\<f-P_0,f-P_0\>}
\<f-P_0,\xi\>I)+
{\cal P}\nabla_X^\perp \xi}. \,\,\,\qed\end{array}
$$

\noindent {\em Proof of Theorem \ref{prop:pairs}:\/} Define
$$
\tilde{\zeta}=(\bar{\lambda}^2
+\bar{\nu}^2)^{-1/2}(\bar{\nu}{\cal P}\eta-
\bar{\lambda}{\cal P}\zeta)\;\;\;\;\;\mbox{and}\;\;\;\;\;
\tilde{\eta}=(\bar{\lambda}^2
+\bar{\nu}^2)^{-1/2}(\bar{\lambda}{\cal P}\eta+\bar{\nu}{\cal P}\zeta)
$$
where  $R^2\bar{\lambda}= \lambda\|\phi\|^2+2\<\phi,\eta\>$, $R^2\bar{\nu}= 2\<\phi,\zeta\>$
and
${\cal P}$ is given by (\ref{eq:vbi}).
Using (\ref{eq:sop}),  we obtain that the shape operators $\tilde{A}_{\tilde{\eta}}$ and
$\tilde{A}_{\tilde{\zeta}}$ of
$
\tilde{\phi}= {\cal I}\circ \phi
=\frac{R^2}{\|\phi\|^2}\phi
$
are given as in (\ref{sffa}) with $\lambda$ and $\mu$ replaced, respectively, by
$$
{\displaystyle \tilde{\lambda}=(\bar{\lambda}^2
+\bar{\nu}^2)^{1/2}\;\;\;\;\mbox{and}\;\;\;\;\tilde{\mu}=\frac{1}{R^2}\|\phi\|^2\mu.}
$$
The pair $(\tilde{g},\tilde{h})$ of conjugate minimal surfaces associated to $\tilde{\phi}$ is
\be\label{eq:til}
\tilde{g}=\tilde{\phi}+\tilde{r}\tilde{\eta}
=\tilde{\phi}+\frac{\bar{\lambda}{\cal P}\eta
+\bar{\nu}{\cal P}\zeta}{\bar{\lambda}^2+\bar{\nu}^2}
\;\;\;\;\mbox{and}\;\;\;\;
\tilde{h}=-\tilde{r}\tilde{\zeta}
=-\frac{\bar{\nu}{\cal P}\eta
-\bar{\lambda}{\cal P}\zeta}{\bar{\lambda}^2
+\bar{\nu}^2}.
\ee
We have
\be\label{eq:r4}
R^4(\bar{\nu}^2+\bar{\lambda}^2)=4(\<\phi,\zeta\>^2
+\<\phi,\eta\>^2+\lambda\<\phi,\eta\>\|\phi\|^2)
+\lambda^2\|\phi\|^4.
\ee
On the other hand, from
$$
{\cal P}\eta
=\eta-2\frac{\<\phi,\eta\>}{\|\phi\|^2}\phi\;\;\;\;
\mbox{and}\;\;\;\;
{\cal P}\zeta
=\zeta-2\frac{\<\phi,\zeta\>}{\|\phi\|^2}\phi,
$$
we obtain
\be\label{eq:peta}
 R^2(\bar{\lambda}{\cal P}\eta
+\bar{\nu}{\cal P}\zeta)
=2\<\phi,\zeta\>\zeta
+2\<\phi,\eta\>(\eta-\lambda \phi)
+\lambda\|\phi\|^2\eta
-\frac{4(\<\phi,\eta\>^2+\<\phi,\zeta\>^2)}{\|\phi\|^2}\phi
\ee
and
\be\label{eq:pzeta}
R^2(\bar{\nu}{\cal P}\eta-\bar{\lambda}{\cal P}\zeta)
=2\<\phi,\zeta\>(\eta+\lambda \phi)
-(2\<\phi,\eta\>+\lambda\|\phi\|^2)\zeta.
\ee
Using that $\lambda=1/r$, $\phi=g-r\eta$ and $h=-r\zeta$,
we have
$$
\<\phi,\zeta\>=-\frac{1}{r}\<g,h\>,\;\;\;\;
\<\phi,\eta\>=\<g,\eta\>-r
\;\;\;\mbox{and}\;\;\;\|\phi\|^2=\|g\|^2-2r\<g,\eta\>+r^2.
$$
Thus,
$$
\bar{\lambda}^2+\bar{\nu}^2=\frac{1}{R^4\|h\|^2}(4\<g,h\>^2
+(\|g\|^2-\|h\|^2)^2).
$$
Hence, from (\ref{eq:til}) we get
$$
\frac{1}{R^2}\tilde{g}= \frac{\phi}{\|\phi\|^2}
+\frac{\bar{\lambda}{\cal P}\eta+\bar{\nu}{\cal P}\zeta}{R^2(\bar{\nu}^2+\bar{\lambda}^2)}
=\frac{(\|g\|^2-\|h\|^2)g+2\<g,h\>h}{4\<g,h\>^2+(\|g\|^2-\|h\|^2)^2}
$$
and
$$
\frac{1}{R^2}\tilde{h}= \frac{2\<g,h\>g
-(\|g\|^2-\|h\|^2)h}{4\<g,h\>^2
+(\|g\|^2-\|h\|^2)^2}.
$$
Therefore,
$$
(\tilde{g},\tilde{h})=F\circ (g,h),
$$
where $F\colon\,\R^{n+2}\times \R^{n+2}\to \R^{n+2}\times \R^{n+2}$ is given by
$$
\frac{1}{R^2}F(U,V)
=\frac{((\|U\|^2-\|V\|^2)U+2\<U,V\>V,2\<U,V\>U
-(\|U\|^2-\|V\|^2)V)}
{4\<U,V\>^2+(\|U\|^2-\|V\|^2)^2}.
$$
Then, as a map $F\colon\,\C^{n+2}\to \C^{n+2}$,
we can write $F$  as
$$
\frac{1}{R^2}F(Z)=\frac{\<\!\< Z,Z\>\!\>}{|\<\!\< Z,Z\>\!\>|^{\,2}}\bar{Z}=\frac{\bar{Z}}{\overline{\<\!\< Z,Z\>\!\>}}=
\overline{T(Z)}.\qed
$$
\begin{remark}\po\label{pairs} {\em For an inversion ${\cal I}$ with respect to a sphere
centered at an arbitrary
point $P_0\in \R^{4}$, the formula in Theorem \ref{prop:pairs} becomes
$\tilde{g}-P_0-i\tilde{h}={T}_R\circ (g-P_0+ih).$
}
\end{remark}

\section[Superconformal surfaces and holomorphic curves]{Superconformal surfaces
and holomorphic curves}

The main goal of this section is to prove Theorems \ref{cor:surf2} and \ref{duality}.
Along the way we establish some additional facts on holomorphic curves and their compositions
with an inversion in $\R^4$. We also look at the Whitney sphere in the light of our results.

\subsection[A characterization of holomorphic curves]{A characterization of holomorphic curves}

In this subsection we prove the following characterization of holomorphic curves that
is interesting on its own right.

\begin{proposition}\po\label{le:hol}
Let $g\colon M^2\to \R^4$ be a simply connected oriented minimal
surface whose holomorphic representative $g+ih\colon\,M^2\to\C^4$
takes values in the quadric $\Q_0\subset\C^4$. Then, $g$ is
holomorphic with respect to some complex structure ${\cal J}$ on
$\R^4$ and its conjugate minimal surface is $h={\cal J}g$.
\end{proposition}
\proof   By the assumption we have
\be\label{cond}
\<g,h\>=0\;\;\;\;\;\mbox{and}\;\;\;\;\;\|g\|^2=\|h\|^2.
\ee
Differentiating the first of  equations (\ref{cond}) and using
that  $h_*=g_*\circ J$, we obtain that the tangent
components $g^T$ and $h^T$ of the position vectors of $g$ and $h$,
respectively, are related by $h^T=Jg^T$. Then, by (\ref{cond}) there
exists a complex structure $\hat{J}$ on $T_g^\perp M$ such that
$\hat{J}g^N=h^N$.

   Differentiating (\ref{cond}) twice gives
$$
\<\alpha_h(X,Y),g^N\>=-\<\alpha_g(X,Y),h^N\>
=-\<\alpha_g(X,Y),\hat{J}g^N\>
=\<\hat{J}\alpha_g(X,Y),g^N\>
$$
and
$$
\<\alpha_h(X,Y),h^N\>=\<\alpha_g(X,Y),g^N\>=
-\<\alpha_g(X,Y),\hat{J}h^N\>
=\<\hat{J}\alpha_g(X,Y),h^N\>.
$$
Hence, $\alpha_h=\hat{J}\circ \alpha_g$.
Since conjugate minimal surfaces satisfy
$\alpha_h(X,Y)=\alpha_g(X,JY)$, it follows that
$$
\alpha_g(X,JY)=\hat{J}\alpha_g(X,Y).
$$
Thus, the complex structure ${\cal J}=J\oplus \hat{J}$ on $g^*T\R^4$
extends to a  complex structure  on $\R^4$  with respect to
which  $g$ is
holomorphic,  and  we have that $h={\cal J}g$.\qed

\subsection{Minimal surfaces and
inversions of holomorphic curves}

We now determine the holomorphic representative of the minimal
surface associated to  the composition of an
inversion with a holomorphic curve. Given an oriented
holomorphic curve \mbox{$f\colon\,M^2\to \R^{4}$}, we denote by
$\hat{J}$ the complex structure on $T_f^\perp M$ determined by the
opposite orientation to that induced by the vector bundle isometry ${\cal
P}\colon\,T_f^\perp M\to  T_{{\cal I}\circ f}^\perp M$ from the
orientation on $T_{{\cal I}\circ f}^\perp M$ defined as in the proof of
Theorem \ref{thm:surf}.

\begin{proposition}\po\label{prop:hol1} Let $f\colon\,M^2\to \R^{4}$ be
a holomorphic curve and let ${\cal I}$ be the inversion in $\R^{4}$ with
respect to the sphere of radius $R$ centered at the origin. Then, the
holomorphic curve in $\C^4$
associated to $\tilde{f}={\cal I}\circ f$ is
\be\label{eq:hol}
\tilde{g}+i\tilde{h}=\frac{R^2}{2\|f^N\|^2}(f^N+i\hat{J}f^N).
\ee
\end{proposition}

\proof Following the proof of Theorem  \ref{prop:pairs} with
$\lambda=0$, we obtain from (\ref{eq:r4}) that
$$
R^4(\bar{\nu}^2+\bar{\lambda}^2)=4\|f^N\|^2.
$$
On the other hand,  (\ref{eq:peta}) and (\ref{eq:pzeta}) now reduce,
respectively, to
\be\label{eq:P1}
{\displaystyle R^2(\bar{\lambda}{\cal P}\eta
+\bar{\nu}{\cal P}\zeta)}
={\displaystyle 2f^N-4\frac{\|f^N\|^2}{\|f\|^2}f}
\ee
and
\be\label{eq:P2}
R^2(\bar{\nu}{\cal P}\eta-\bar{\lambda}{\cal
P}\zeta))=2{\<f,\zeta\>}\eta-2{\<f,\eta\>}\zeta\\
=-2\hat{J}f^N.
\ee
  Then (\ref{eq:hol}) follows from (\ref{eq:til}), (\ref{eq:P1}) and
(\ref{eq:P2}).\qed

\begin{remark}\po\label{pairs2} {\em For an inversion ${\cal I}$ with
respect to a sphere centered at an
arbitrary point $P_0\in \R^{4}$, formula (\ref{eq:hol}) becomes
$$
\tilde{g}+i\tilde{h}=
P_0+\frac{R^2}{2\|(f-P_0)^N\|^2}((f-P_0)^N + i\hat{J}(f-P_0)^N),
$$
where  $(f-P_0)^N$ is the normal component of the position vector of $f$ with respect to $P_0$.}
\end{remark}

The first assertion in Theorem \ref{duality} is a consequence of the following general fact.

\begin{proposition}\po\label{prop:Kaehler} Let $f\colon\,M^n\to\C^{n+p}$ be a holomorphic
isometric immersion of
a Kaehler manifold of real dimension $2n$. Then
$F=f^N/\|f^N\|^2$ is anti-holomorphic with respect to the complex structure $J$ of $\C^{n+p}$
if and only if $f^T\in\kerl\alpha_f^\perp$, where $\alpha_f^\perp$ is the
$\spa\{f^N,Jf^N\}^\perp$ component
of $\alpha_f$. In particular, this is always the case if the complex codimension is $p=1$.
\end{proposition}

\proof Denote for simplicity  $g=f^N$ and $h=f^T$.  We have
$$
g_*Z=Z-\nabla_Zh-\alpha_f(Z,h)
$$
and
$$
g_*Z=-A^f_gZ+\nabla^\perp_Zg.
$$
We conclude that
$$
g_*Z=-A^f_gZ-\alpha_f(Z,h).
$$
Then
$$
F_*Z=\frac{1}{\|g\|^4}(\|g\|^2g_*Z-2\<g_*Z,g\>g),
$$
and hence
$$
\|g\|^4F_*Z=-\|g\|^2(A^f_gZ+\alpha_f(Z,h))
+2\<\alpha_f(Z,h),g\>g.
$$
Equivalently,
$$
\|g\|^2F_*Z=-A^f_gZ+\frac{1}{\|g\|^2}
(\<\alpha_f(Z,h),g\>g-\<\alpha_f(Z,h),Jg\>Jg)
-\alpha^\perp_f(Z,h).
$$
It is now easy to see that
$$
F_*\circ J=-J\circ F_*\;\;\;\;\mbox{if and only if}\;\;\;\;
h\in\kerl\alpha_f^\perp.\,\,\,\,\qed
$$

\begin{remark}\po\label{hatJ} {\em Now that we know that for a holomorphic
curve $f\colon\,M^2\to \R^4$ with respect to a complex structure $J$ on $\R^4$
the map $f^N/\|f^N\|^2$ defines an anti-holomorphic curve with respect to $J$,
it follows that the complex structure $\hat{J}$ on $T_f^\perp M$ is the
restriction of $-J$ to $T_f^\perp M$.}
\end{remark}

\subsection[Proof of Theorems \ref{cor:surf2} and \ref{duality}]
{Proof of Theorems \ref{cor:surf2} and \ref{duality}}
\noindent {\em Proof of  Theorem \ref{cor:surf2}:\/}
 If either $\phi$ or
its dual  parameterizes the composition
of an inversion with respect to a sphere of radius $R$
centered, say, at the origin, with a holomorphic curve with respect to a complex structure $J$, then the associated minimal surface is $g=\frac{R^2f^N}{2\|f^N\|^2}$, which is anti-holomorphic with respect to $J$ by the assertion in Theorem \ref{duality} already proved. If either $\phi$
or its dual only takes the value  $0$,
then the holomorphic representative $g+ih$ of the associated minimal
surface $g$ takes values in $\Q_0$, and hence $(iii)$ follows from
 Proposition~\ref{le:hol}.
Conversely, if $g$ is  holomorphic with respect to some complex structure
${\cal J}$ on $\R^4$, $h={\cal J}g$   and $\hat{J}_+$ is the  restriction of $-{\cal J}$
to $T_g^\perp M$, then $\phi_-$ reduces to the constant map $\phi_-=0$
and its dual is $\phi_+=2g^N$. Thus,  the assertion in Theorem \ref{duality}
already proved implies that $\phi=\phi_+=2g^N$ is
 as in $(i)$.\vspace{2ex}\qed

\noindent {\em Proof of  Theorem \ref{duality}:\/} We already know that
$f\mapsto f^*=f^N/2\|f^N\|^2$
maps  ${\cal H}_+$ into ${\cal H}_-$. Since $f^*$ is the minimal surface
associated to $f/\|f\|^2$ by Proposition \ref{prop:hol1},
that $f$ and $f^*$ induce conformal surfaces on $L^2$ follows
from Theorem~\ref{thm:surf}. Finally, since $f^*$ is anti-holomorphic with respect to $J$,
the fact that  $f^*$ is the minimal surface associated to $f/\|f\|^2$ implies that
$f/\|f\|^2=2(f^*)^N$.
Hence,
$$
f=\frac{2(f^*)^N}{\|2(f^*)^N\|^2}
=\frac{(f^*)^N}{2\|(f^*)^N\|^2}=(f^*)^*.
$$
In particular, this implies that the map $f\mapsto f^*$ is a bijection between
${\cal H}_+$ and ${\cal H}_-$.\qed

\subsection[The  Whitney sphere]{The  Whitney sphere}

The {\em Whitney sphere\/} is the immersion $$
(x,y,z)\in\Sf^2\mapsto\frac{1}{1 + z^2}(x(1+iz),y(1+iz))
\in \C^2
$$
of  the unit sphere $\Sf^2(1)$ into the complex Euclidean plane. Among
several interesting characterizations,  the
one that concerns this paper is as the only Lagrangean superconformal
surface in
$\C^2$ due to Castro \cite{ca}.

The Whitney sphere is just the composition ${\cal I}\circ f$ of the
holomorphic curve $f\colon\,\C^*\to \C^2$ given by
$$
f(z)=\left(z,1/z\right)
$$
with the inversion ${\cal I}$ with respect to the sphere of unit radius centered at
the origin.

By Proposition \ref{prop:hol1}, if we orient $\Sf^2$ and $T_{{\cal
I}\circ f}^\perp \Sf^2$ as described in the
paragraph preceding it, then the pair $(g,h)$ of conjugate minimal
surfaces associated to ${\cal I}\circ f$ is given by
$$
g=\frac{1}{2}{\cal I}\circ f^N\,\,\,\,\mbox{and}\,\,\,\,
h=\frac{i}{2}{\cal I}\circ f^N.
$$
An easy computation shows that
$$
g(z)=\frac{1}{4}\left(1/\bar{z},\bar{z}\right)\,\,\,\,\mbox{and}\,\,\,\,h(z)=
\frac{i}{4}\left(1/\bar{z},\bar{z}\right).
$$
In fact, the construction of Theorem \ref{thm:surf} applied to $g$ gives
$$
g+{\cal J}_-h=2g^N=\frac{1}{2}{\cal I}\circ
\frac{1}{2}f={\cal I}\circ f.
$$

\section[Superconformal and superminimal surfaces ]{Superconformal and superminimal surfaces}

In this  section we prove Theorem \ref{cor:surf3} and
illustrate our result with the Veronese surface.

\subsection[Proof of Theorem \ref{cor:surf3}]{Proof of Theorem \ref{cor:surf3}}

First we consider the case of $\Sf_R^4=\Sf^4(Re_5;R)$. Regard the
stereographic projection of $\Sf_R^4$  onto $\R^4$ as the
restriction to $\Sf_R^4$ of the inversion ${\cal I}$ in $\R^5$ with respect to the
sphere of
radius $2R$ centered at $2Re_5$. Now, given a superminimal surface
\mbox{$f\colon\,M^2\to \Sf_R^4$}, we obtain as in Proposition
\ref{prop:hol1} (see also Remark~\ref{pairs2}) that the pair
$(g,h)$ of conjugate minimal surfaces associated to
${\cal I}\circ f$ is
$$
{g}=2Re_5+\frac{2R^2(f-2Re_5)^N}{\|(f-2Re_5)^N\|^2}\,\,\,\,\,\,\mbox{and}\,\,\,\,\,\,\,{h}
=\frac{2R^2\hat{J}(f-2Re_5)^N}{\|(f-2Re_5)^N\|^2},
$$
where  $(f-2Re_5)^N$ is the normal component of $f$ with respect to $2Re_5$. Thus,
\be\label{eq:gh}\<g,h\>=0\,\,\,\,\,\mbox{and}\,\,\,\,\,\|g-2Re_5\|=\|h\|,\ee
and hence $g+ih$ takes values in $\Q_{4R^2}$.

     Conversely, assume that the pair   $(g,h)$ of conjugate minimal surfaces that
     gives rise to $\phi_+$ and $\phi_-$
satisfies (\ref{eq:gh}). We claim that ${\cal I}\circ \phi_\pm$ is a superminimal surface
in $\Sf_R^4$.
The arguments for  $\phi_+$ and $\phi_-$ being the same, in the sequel we omit the
subscript $\pm$ for simplicity of notation.

 Set $\bar{R}=2R$, $P_0=2Re_5$ and define
$$
\tilde{\zeta}=\frac{\bar{\nu}{\cal P}\eta- \bar{\lambda}{\cal
P}\zeta}{(\bar{\lambda}^2 +\bar{\nu}^2)^{1/2}}
\;\;\;\;\mbox{and}\;\;\;\; \tilde{\eta}=\frac{\bar{\lambda}{\cal
P}\eta +\bar{\nu}{\cal P}\zeta}{(\bar{\lambda}^2
+\bar{\nu}^2)^{1/2}},
$$
where ${\cal P}$ is the vector bundle isometry between
$T_f^\perp M$ and $T_{{\cal I}\circ f}^\perp M$ given by (\ref{eq:vbi}),
$$
\bar{R}^2\bar{\nu}=
2\<\phi-P_0,\zeta\>\;\;\;\;\;\mbox{and}\;\;\;\;\;
\bar{R}^2\bar{\lambda}= \lambda\|\phi-P_0\|^2+ 2\<\phi-P_0,\eta\>.
$$
Using Lemma \ref{le:inv}, we obtain that the shape operators
$\tilde{A}_{\tilde{\eta}}$ and $\tilde{A}_{\tilde{\zeta}}$ of $
{\cal I}\circ \phi$ are given as in (\ref{sffa}) with $\lambda$
and $\mu$ replaced, respectively, by
$$
\tilde{\lambda}=(\bar{\lambda}^2+\bar{\nu}^2)^{1/2}\,\,\,\,\,\,\mbox{and}\,\,\,\,\,
\,\tilde{\mu}=\frac{\|\phi-P_0\|^2}{\bar{R}^2}\mu.
$$
Since  (\ref{eq:gh}) holds, using that
${\displaystyle h=-\frac{1}{\lambda}\zeta}$ and ${\displaystyle
g-P_0=\phi-P_0+\frac{1}{\lambda}\eta}$  we obtain
$$
\<\phi-P_0,\zeta\>=\<g-P_0-\frac{1}{\lambda}\eta,\zeta\>=0$$ and
$$
\frac{-2}{\lambda}\<\phi-P_0,\eta\>=\|\phi-P_0\|^2+\frac{1}{\lambda^2}-\|g-P_0\|^2=\|\phi-P_0\|^2+
\frac{1}{\lambda^2}-\|h\|^2=\|\phi-P_0\|^2.$$ Thus,
$\bar{\nu}=0=\bar{\lambda}$, and hence $\tilde{\lambda}=0$, which
is equivalent to ${\cal I}\circ \phi$ being superminimal.

   In the case of $\Hy_R^4=\Hy^4(-Re_5;R)$, let ${\cal I}\colon\,\Hy^4_R\to \R^4$
   denote the
 stereographic projection defined in the introduction. Notice that
    the image of $\Hy^4_R$ by ${\cal I}$ is the
   open ball $B(0;2R)\subset \R^4$,  the   induced metric on $B(0;2R)$ being
   the Poincar\'e hyperbolic metric of constant sectional curvature $-1/R^2$.
Observe also  that ${\cal I}$ can  be regarded as the
   restriction to  $\Hy^4_R$ of the ``inversion" ${\cal I}$ on $\Les^5$ with
   respect to $\Hy^4(-2Re_5;2R)$ (see Lemma \ref{le:invv}).
Taking Lemma \ref{le:invv} into account, the remaining of the proof is entirely
similar to that of the spherical case, the conclusion now being  that $(g,h)$
is the pair of conjugate minimal surfaces associated to ${\cal I}\circ
f$, where
$f\colon\,M^2\to \Hy^4_R$ is a superminimal surface, if and only if
\be\label{eq:gh2}\<g,h\>=0\,\,\,\,\,\mbox{and}\,\,\,\,\,\<g+2Re_5,g+2Re_5\>=\<h,h\>,\ee
that is, if and only if  $g+ih$ takes values in $\Q_{-4R^2}$. \qed

\subsection[The  Veronese surface]{The  Veronese surface}

The Veronese surface
$f\colon\,\Sf^2_{1/3}\to\Sf^4\subset\R^5$ given by
$$
f(x,y,z)=\frac{1}{2\sqrt{3}}
(2xy,2xz,2yz,x^2-y^2,\frac{1}{\sqrt{3}}(x^2+y^2-2z^2))
$$
is a superminimal surface with constant normal curvature in the sphere.
In  spherical coordinates
$(x,y,z)=\sqrt{3}(\sin\va\cos\theta,\sin\va\sin\theta,
\cos\va)$, we can write $f$ as
$$
\frac{2}{\sqrt{3}}f(\va,\theta)
=\sin^2\va X_1+\sin 2\va X_2 +\frac{1}{\sqrt{3}}(1-3\cos^2\va)e_5,
$$
where $(\va,\theta)\in (0,\pi)\times (0,2\pi)$ and
$$
X_1=\sin 2\theta e_1+\cos 2\theta e_4,\;\;\;\;\; X_2=\cos\theta
e_2+\sin\theta e_3.
$$
We have
$$
\frac{2}{\sqrt{3}}\frac{\d f}{\d
\va}=\sin 2\va X_1+2\cos 2\va X_2+{\sqrt{3}}\sin 2\va e_5
$$
and
$$
\frac{2}{\sqrt{3}}\frac{\d f}{\d \theta}=2\sin^2\va X_3+\sin 2\va X_4,
$$
where
$$
X_3=\cos 2\theta e_1-\sin 2\theta e_4,\;\;\;\;\;
X_4=-\sin\theta e_2+\cos\theta e_3.
$$
Thus, an orthonormal frame $\{\eta,\zeta\}$ of $T_f^\perp \Sf^2_{1/3}$
is given by
$$
\eta=-\cos\va X_3+\sin\va X_4,\;\;\;\;\;
2\zeta=(1+\cos^2\va)X_1-\sin 2\va X_2
-\sqrt{3}\sin^2\va e_5.
$$
Identifying $\R^4$ with the tangent
space of $\Sf(e_5,1)\subset \R^5$ at the origin and viewing $f$ as a map into
$\Sf(e_5,1)$, the pair $(g,h)$ of conjugate minimal surfaces associated to the
stereographic projection of $f$ onto $\R^4$ is
$$g=2e_5+2\frac{(f-2e_5)^N}{\|(f-2e_5)^N\|^2}\;\;\;\;\;\mbox{and}\;\;\;\;\;h=2\frac{\hat{J}(f-2e_5)^N}{\|(f-2e_5)^N\|^2},$$
where $(f-2e_5)^N$ denotes the component of the position vector
$f-2e_5$ in the normal bundle of $f$ (as a map into $\Sf(e_5,1)$).
Using that $(f-e_5)^N=0$, we obtain
$$
\frac{1}{2}g=e_5-\frac{\<\eta,e_5\>\eta+\<\zeta,e_5\>\zeta}{\<\eta,e_5\>^2+\<\zeta,e_5\>^2}=
e_5-\frac{\zeta}{\<\zeta,e_5\>}\;\;\;\;\;\mbox{and}\;\;\;\;\;\frac{1}{2}h=-\frac{\eta}{\<\zeta,e_5\>}.
$$
Therefore, we have the conjugate minimal surfaces
$$
g=\frac{2}{\sqrt{3}\sin^2\va}\left((1+\cos^2\va)X_1
-2\sin\va\cos\va X_2\right)\;\;\mbox{and}\;\;
h=\frac{4}{\sqrt{3}\sin^2\va}\left(\cos\va X_3-\sin \va
X_4\right),
$$
which induce the complete metric
$$
ds^2=\frac{4(1+3\cos^2\va)}{\sin^6\va}(\sin^2\va
d^2\theta+d^2\va).
$$
Amazingly enough, in this case we have a compact superconformal surface generated
as in (\ref{eq:phi3}) by a complete unbounded minimal surface.

\vspace{.5in} {\renewcommand{\baselinestretch}{1}

\hspace*{-20ex}\begin{tabbing} \indent\= IMPA -- Estrada Dona Castorina, 110
\indent\indent\=  Universidade Federal de S\~{a}o Carlos\\
\> 22460-320 -- Rio de Janeiro -- Brazil  \>
13565-905 -- S\~{a}o Carlos -- Brazil \\
\> E-mail: marcos@impa.br \> E-mail: tojeiro@dm.ufscar.br
\end{tabbing}}

\end{document}